\documentclass[letterpaper, 12pt]{article} 

\usepackage[paper=letterpaper,left=19mm,right=19mm,top=25.4mm,bottom=19mm,footskip=5.3mm]{geometry}
\usepackage{setspace, amsmath, amssymb, graphicx, framed, mathrsfs, rotating, color, multirow, amsthm}
\usepackage{wrapfig, subfigure, enumerate, caption}
\usepackage[numbers,sort&compress]{natbib}
\usepackage{mathptmx, titlesec}  
\usepackage[absolute]{textpos} 

\usepackage[pdftex]{hyperref}

\hypersetup{
pdftitle={Parameter Optimization of a Virtual Synchronous Machine in a Microgrid},
pdfsubject={Article},
pdfauthor={Timo Dewenter, Wiebke Heins, Benjamin Werther, Alexander K. Hartmann, Christian Bohn, Hans-Peter Beck},
pdfkeywords={Optimization, MicroGrid},
colorlinks=true,
urlcolor=black,
linkcolor=black,
citecolor=black
}

\makeatletter
\g@addto@macro\normalsize{%
\setlength{\abovedisplayshortskip}{0.2\baselineskip}
\setlength{\belowdisplayshortskip}{0.2\baselineskip}
\setlength{\abovedisplayskip}{0.2\baselineskip}
\setlength{\belowdisplayskip}{0.2\baselineskip}}
\makeatother

\makeatletter
\newcommand{\MSonehalfspacing}{%
	\setstretch{1.44}
	\ifcase \@ptsize \relax 
		\setstretch {1.448}%
	\or 
		\setstretch {1.399}%
	\or 
		\setstretch {1.433}%
	\fi
}
\newcommand{\MSdoublespacing}{%
	\setstretch {1.92}
	\ifcase \@ptsize \relax 
		\setstretch {1.936}%
	\or 
		\setstretch {1.866}%
	\or 
		\setstretch {1.902}%
	\fi
}

\renewcommand\section{\@startsection
	{section}{1}{0mm}
	{0.4\baselineskip}
	{0.01\baselineskip}
	{\normalfont\bfseries\normalsize}
}
\makeatother
\makeatletter
\renewcommand\subsection{\@startsection
	{subsection}{2}{0mm}
	{0.4\baselineskip}
	{0.01\baselineskip}
	{\normalfont\bfseries\normalsize}
}
\makeatother
\makeatletter
\renewcommand\subsubsection{\@startsection
	{subsubsection}{3}{0mm}
	{0.2\baselineskip}
	{0.2\baselineskip}
	{\normalfont\itshape\normalsize}
}

\makeatother

\MSdoublespacing

\newcommand{\mathsym}[1]{{}}

\newcommand{\FG}[1]{Fig.~\ref{#1}}
\newcommand{\mytitlesize}{\fontsize{20pt}{20pt}\selectfont}
\newcommand{\myfootnotesize}{\fontsize{9pt}{9pt}\selectfont}
\newcommand{\T}[1]{\text{#1}}

\begin{document}

\begin{textblock*}{\textwidth}(19mm,1.1\textheight)
\noindent \myfootnotesize \textit{Manuscript received 23 June 2016}
\end{textblock*}

\begin{center}
\textbf{\mytitlesize{PARAMETER OPTIMISATION OF A VIRTUAL SYNCHRONOUS MACHINE IN A MICROGRID}}
\end{center}
\begin{center}
\singlespacing
{\small Timo Dewenter,$^1$ Wiebke Heins,$^{2,3}$ Benjamin Werther,$^4$ Alexander K. Hartmann,$^1$ Christian Bohn,$^2$ \\ and Hans-Peter Beck$^4$} \\[1mm]
{\small
\textit{$^1$ Institut f\"ur Physik, Carl von Ossietzky Universit\"at Oldenburg, D-26111 Oldenburg, Germany} \\
\textit{$^2$ Institut f\"ur Elektrische Informationstechnik, TU Clausthal, D-38678 Clausthal-Zellerfeld, Germany} \\
\textit{$^3$ Zentrum f\"ur Technomathematik, Universit\"at Bremen, D-28359 Bremen, Germany} \\
\textit{$^4$ Institut f\"ur Elektrische Energietechnik und Energiesysteme, TU Clausthal, \\ D-38678 Clausthal-Zellerfeld, Germany}}
\end{center}

\section*{Abstract}
\noindent Parameters of a virtual synchronous machine in a small microgrid are optimised. The dynamical behaviour of the system is simulated after a perturbation, where the system needs to return to its steady state. The cost functional evaluates the system behaviour for different parameters. This functional is minimised by Parallel Tempering. 
Two perturbation scenarios are investigated and the resulting optimal parameters agree with analytical predictions. Dependent on the focus of the optimisation different optima are obtained for each perturbation scenario. During the transient the system leaves the allowed voltage and frequency bands only for a short time if the perturbation is within a certain range.
\thispagestyle{empty}
\setcounter{page}{1}

\section*{Key Words}
\noindent Inverter-Based Microgrid, Virtual Synchronous Machine, Stochastic Optimisation, Parallel Tempering 

\section{Introduction}
\noindent The number of renewable distributed energy sources (DER) has increased in the last decades forced by political, ecological, and economical aspects. Many DER are attached to the low-voltage grid by inverters, whose increased usage is accompanied by the need to find suitable control strategies and parameters for, e.g., frequency-power droop control in autonomous microgrids.
Simulation methods, models and stability conditions for microgrids based on droop-controlled inverters are investigated in 
\cite{Pogaku2007,Coelho2002,Marwali2007,Barklund2008,Soultanis2007,Engler2005,Simpson_Porco2013}.
A rigorous stability \newpage
\noindent analysis is done in \cite{DroopControlSchiffer2014}, in which conditions on the droop gains are derived.
Simulations \cite{Raghami2015,Godoy2012,Mishra2011,Sanjari2013} have been used to obtain optimal control parameters of inverters or distributed generators in microgrids.
\textit{Particle swarm optimisation} in which a ``swarm'' of solutions moves in the search-space is used in 
\cite{Chung2011,Hassan2011,Al_Saedi2012,Bevrani2012,Yu2014,Zeng2014,Wang2016}.

To enhance stability in microgrids, one can use programmable inverters with storage, as e.g., the virtual synchronous machine (VISMA) \cite{Beck2007}. It is a hysteresis controlled three-phase inverter, whose setpoints are determined by a synchronous machine model implemented on a control computer. Inertia to improve transient stability of the grid, is provided by a storage device. The VISMA is able to control (re-)active power bidirectionally and can be adjusted to meet specific power system requirements.
 
Here, the VISMA as grid-building element in a low-voltage islanded microgrid with voltage source inverters is investigated. The basic control strategy is droop control \cite{Engler2005,DroopControlSchiffer2014} for both voltage and frequency. 
We use the parallel tempering method \cite{Hukushima96,Swendsen86} for optimisation of the VISMA parameters under varying transient loads (see e.g. \cite{RahmanTito2013}).
Parallel Tempering allows to find (near-) optimal solutions for complex optimisation problems \cite{Hartmann01,Hartmann04} efficiently.  
The objective of our analysis is to show that the optimisation method is generally applicable to determine optimal control parameters in microgrids. Furthermore, different types of optima allow insights in the effects of the VISMA in combination with regular droop controlled inverters in microgrids for the first time. 

In Sec.~2, the simulation model is described. The optimisation problem is stated in Sec.~3 and Sec.~4 explains the implementation. Results are presented in Sec.~5, a conclusion is given in Sec.~6.

\section{Model of an Inverter-Based Microgrid with VISMA}
\subsection{Lines and Loads}
\noindent Lines are modelled as algebraic equations describing the relation between voltage angles $\theta_i(t)$ and voltage magnitudes $V_i(t)$ at grid node $i \in [1,n]$ and (re-)active power flows \cite{Kundur1994}. Magnitudes and angles for all grid nodes are gathered in $\boldsymbol{V}(t) = \left[ V_1(t),~V_2(t),~...~,~V_n(t) \right]^{\textnormal{T}}$ and $\boldsymbol{\theta}(t) = \left[ \theta_1(t),~\theta_2(t),~...~,~\theta_n(t) \right]^{\textnormal{T}}$.
The (re-)active power injected at node $i$ is then described by the power balance equations:
\begin{align}
P_i (\boldsymbol{V}(t), \boldsymbol{\theta}(t)) &=&
3 \left( G_{ii} V_i^2(t) - \sum_{k \in N(i)} 
V_i(t) V_k(t) 
\left( G_{ik} \cos(\theta_i(t) - \theta_k(t)) 
+ B_{ik} \sin(\theta_i(t) - \theta_k(t)) \right) \right) \\
Q_i (\boldsymbol{V}(t), \boldsymbol{\theta}(t)) &=&
3 \left( -B_{ii} V_i^2(t) - \sum_{k \in N(i)} 
V_i(t) V_k(t) 
\left( G_{ik} \sin(\theta_i(t) - \theta_k(t)) 
- B_{ik} \cos(\theta_i(t) - \theta_k(t)) \right) \right)
\end{align}

Here, $G_{ii} = \hat{G}_{ii} + \sum_{k \in N(i)} G_{ik}$ and $B_{ii} = \hat{B}_{ii} + \sum_{k \in N(i)} B_{ik}$, and $k \in N(i)$ denotes summation over neighbours~$k$ of node~$i$, $\underline{\hat{Y}}_{ii}=\hat{G}_{ii}+\textnormal{j} \hat{B}_{ii}$ its shunt admittance, and $\underline{Y}_{ik}=G_{ik}+\textnormal{j} B_{ik}$ the admittance of line $ik$.

The load is modelled as an external disturbance $S_{\textnormal{load}}(t) = P_{\textnormal{load}}(t) + \textnormal{j} Q_{\textnormal{load}}(t)$.
An algebraic constraint is introduced for the node $k$ to which the load is connected, so 
$P_{\textnormal{load}}(t) = P_k(\boldsymbol{V}, \boldsymbol{\theta})$, $Q_{\textnormal{load}}(t) = Q_k(\boldsymbol{V}, \boldsymbol{\theta})$.

\subsection{Droop-Controlled Inverters \label{sec:droop-controlled-inverters}}
\noindent Following \cite{DroopControlSchiffer2014}, inverters are modelled as controllable AC voltage sources described by differential equations.
Each inverter is connected to the grid via an LCL-filter with inductance $L_{\textnormal{inv}}$ on the inverter side, filter capacitance $C_{\textnormal{f}}$ and coupling inductance $L_{\textnormal{C}}$.
Here, $V_i(t)$ and $\theta_i(t)$ denote time-varying voltage magnitudes and angles over filter capacitances $C_{\textnormal{f}}$, assuming that these are the voltages controlled by the inverter.

Droop frequency and voltage control is based on decentralized proportional controllers. Its adaption to inverter-based microgrids has been investigated 
\cite{Engler2005, Barklund2008, Chandorkar1993, Coelho2002, Soultanis2007} extensively. Because droop control is purely proportional, an offset error occurs as soon as the system is  permanently disturbed. The objective of the control strategy is that in steady state (denoted by $\ast$) of the closed loop system devices participating in droop control share the additional (re-)active power caused by the disturbance according to the equations:  
\begin{equation}
k_{\textnormal{P},i} \left( P_{\textnormal{nom},i}- P_i^{\ast}(\boldsymbol{V}^{\ast}, \boldsymbol{\theta}^{\ast}) \right) = \omega_i^{\ast} - \omega_{\textnormal{nom}}, \hspace{0.5cm}
k_{\textnormal{Q},i} \left( Q_{\textnormal{nom},i}- Q_i^{\ast}(\boldsymbol{V}^{\ast}, \boldsymbol{\theta}^{\ast})  \right) = V_i^{\ast} - V_{\textnormal{nom}} \label{eq:droopPQ}
\end{equation} 

Here, $P_{\textnormal{nom},i}$ and $Q_{\textnormal{nom},i}$ denote the nominal active and reactive power injections of each device. $V_{\textnormal{nom}}$ and $\omega_{\textnormal{nom}}$ denote the nominal voltage magnitude and frequency, respectively. The coefficients $k_{\textnormal{P},i}$ and $k_{\textnormal{Q},i}$ are parameters which determine the desired power sharing among devices. 
A common approach for the choice of droop coefficients $k_{\textnormal{P},i}$ and $k_{\textnormal{Q},i}$ is proportional load sharing (see \cite{DroopControlSchiffer2014} for analysis). Based on the power rating $S_i$ of each device and taking into account the legal limits for grid frequency and voltage magnitudes (49.8 Hz to 50.2 Hz, and 207 V to 253 V, respectively), we obtain:
\begin{equation}
k_{\textnormal{P},i} = \frac{0.4 \cdot 2 \pi}{2 S_{i}} \; \frac{\rm rad}{\rm s} = \frac{0.4 \pi}{S_{i}} \; \frac{\rm rad}{\rm s},
\qquad k_{\textnormal{Q},i} = \frac{46 \rm{V}}{2 S_{i}} = \frac{23 \rm{V}}{S_{i}} ~~~ \forall \ i
\label{kp_kQ}
\end{equation}

In \cite{DroopControlSchiffer2014}, voltage source inverters are modelled with instantaneous frequency $\dot{\theta}_i(t) = \omega_{\textnormal{sp},i}(t)$, and first-order-delay voltage control $T_{\textnormal{inv}} \dot{V}_i(t) = - V_i(t) + V_{\textnormal{sp},i}(t)$, where $\omega_{\textnormal{sp},i}(t)$ and $V_{\textnormal{sp},i}(t)$ denote freely adjustable frequency and voltage setpoints. 
Furthermore, power measurements are processed by low-pass filters with time constants $T_i \gg T_{\textnormal{inv}}$. Choosing setpoints $\omega_{\textnormal{sp},i}(t)$ and $V_{\textnormal{sp},i}(t)$ according to \eqref{eq:droopPQ} gives (see \cite{DroopControlSchiffer2014} for details):
\begin{eqnarray}
\dot{\theta}_i(t) &=& \omega_i(t) \label{eq:InverterDTHETA} \\
T_i \dot{\omega}_i(t) &=& 
 - \omega_i(t) + \omega_{\textnormal{nom}}
 + k_{\textnormal{P},i} 
 (P_{\textnormal{nom},i} - P_i(\boldsymbol{V}(t), \boldsymbol{\theta}(t))) \\
T_i \dot{V}_i(t) &=& -V_i(t) + V_{\textnormal{nom}}
 + k_{\textnormal{Q},i} (Q_{\textnormal{nom},i} - Q_i(\boldsymbol{V}(t), \boldsymbol{\theta}(t))) 
\end{eqnarray}

\subsection{The Virtual Synchronous Machine (VISMA) with Droop- and Secondary Frequency Control}
\noindent The VISMA \cite{Beck2007} is a programmable inverter which mimics the dynamics of a synchronous machine.
It uses the three-phase grid voltages as input and the three-phase currents as output.
Its machine model, which defines how the programmable inverter is supposed to act, is adapted to fit in the overall model:
\begin{eqnarray}
  \dot{\theta}_i(t)    &=&   \omega_i(t) \label{eq:dPhi2} \\
  J \dot{\omega}_i(t)  &=&   -\frac{k_{\textnormal{d}}}{T_{\textnormal{d}}} \omega_i(t)
   -\frac{k_{\textnormal{d}}}{T_{\textnormal{d}}} d(t) + \frac{1}{\omega_i(t)} 
   \left( P_{\textnormal{inject}}(t) - P_i(\boldsymbol{V}(t), \boldsymbol{\theta}(t)) \right) \label{eq:VISMAdomega2} \\
   \dot{d}(t) &=& -\frac{1}{T_{\textnormal{d}}} \omega_i(t) -\frac{1}{T_{\textnormal{d}}} d(t) \label{eq:VISMAdamping}
\end{eqnarray}

Parameters are the virtual moment of inertia $J>0$, the mechanical damping factor $k_{\textnormal{d}}>0$, the damping time constant $T_{\textnormal{d}}>0$.
Compared to the VISMA model as stated in \cite{dewenter14}, this model was obtained by defining a "damping state" \mbox{$d(t) = \frac{T_{\textnormal{d}}}{k_{\textnormal{d}}} M_\textnormal{d}(t) - \omega_i(t)$} and replacement of the momentum $M_{\textnormal{mech}}(t)$ by $M_{\textnormal{mech}}(t) = \frac{1}{\omega_i(t)} P_{\textnormal{inject}}(t)$, 
with $P_{\textnormal{inject}}(t)$ denoting the active power injected into the grid by the VISMA. In this setup, it is used for the purpose of droop and secondary frequency control, i.e. $P_{\textnormal{inject}}(t) = P_{\textnormal{droop}}(t) + P_{\textnormal{secondary}}(t)$. According to \eqref{eq:droopPQ} it is $P_{\textnormal{droop}}(t) = P_{\textnormal{nom},i} + \frac{1}{k_{\textnormal{P},i}} (\omega_{\textnormal{nom}} - \omega_i(t))$.
Secondary frequency control is only performed by the VISMA and realized via an integral controller: 
\begin{equation}
\label{eq:integral_controller}
\dot{x}(t) = K_{\textnormal{I}} (\omega_{\textnormal{nom}} - \omega_i(t)), \hspace{0.5cm} P_{\textnormal{secondary}}(t) = x(t) 
\end{equation}

The voltage $E_{\textnormal{P}}$ \cite{dewenter14} is represented here by the voltage magnitude $V_i(t)>0$ of the VISMA. Voltage dynamics of the VISMA are assumed as a first-order delay with a fast time constant $T_{\textnormal{inv},i}$.
A specific voltage control strategy for the VISMA \cite{chen11} is implemented using the root mean square value $V_{\textnormal{grid},i}(t)$ obtained from the grid voltage measurement between stator and grid (cf., \FG{fig:Disturbance_scenario_setup}) as: 
\begin{equation}
T_{\textnormal{inv},i} \dot{V}_i(t) = -V_i(t) + 
V_{\textnormal{nom}} + k_{\textnormal{V}} \left( V_{\textnormal{nom}} - V_{\textnormal{grid},i}(t) \right)
\end{equation}

Furthermore, VISMA stator equations \cite{dewenter14} are simplified as quasi-static and represented via an algebraic equation $\underline{Y}_{\textnormal{VISMA}} = \frac{1}{R_{\textnormal{S}}+\textnormal{j} \omega L_{\textnormal{S}}}$ with stator resistance $R_{\textnormal{S}}>0$ and stator inductance $L_{\textnormal{S}}>0$.
\begin{figure}[ht]
    \centering
    \includegraphics[width=0.85\textwidth]{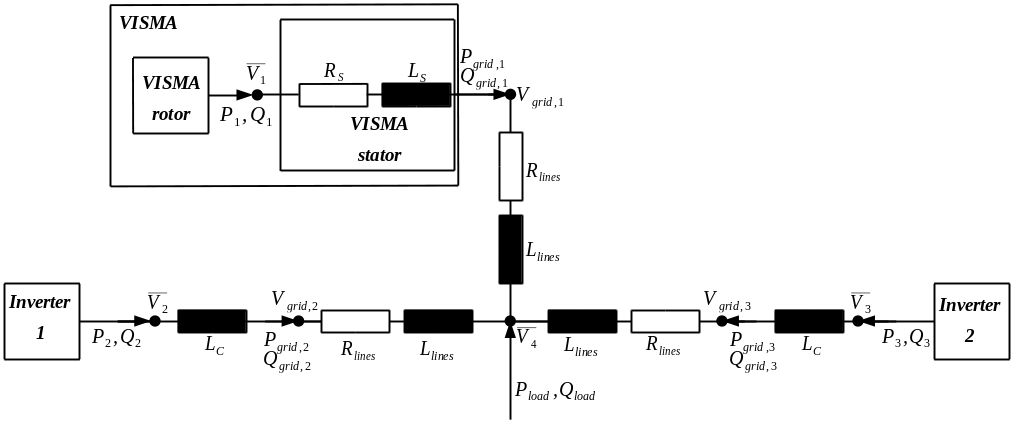}
    \caption{Scheme of the microgrid setup for simulation in perturbation scenarios. \label{fig:Disturbance_scenario_setup}}
\end{figure}

\subsection{Overall Simulation Model with Respect to a Reference Node \label{sec:overall-model-with-respect-to-a-reference-node}}
\noindent We choose the VISMA node (node 1) as reference node. All voltage angles are replaced by their difference to the reference node's voltage angle via $\Delta \theta_i(t) := \theta_i(t) - \theta_1(t) \, \forall \, i$. Naturally, we have $\Delta \theta_1(t) \equiv 0$ and therefore the state $\Delta \theta_1(t)$ and \eqref{eq:dPhi2} are not needed to describe the full system. A new vector is defined for angle states as $\boldsymbol{\Delta \theta}(t) = \left[ \Delta \theta_2(t),~...~,~\Delta \theta_n(t) \right]^{\textnormal{T}} \in \mathbb{R}^{n-1}$. For lines and loads, $\theta_i(t)$ can be directly replaced by $\Delta \theta_i(t) ~ \forall \, i$.
For the inverters, assuming that none of them is connected to node~1, \eqref{eq:InverterDTHETA} is replaced by $\Delta \dot{\theta}_i(t) = \omega_i(t) - \omega_1(t)$.
Given the complex power $\underline{S}_i = P_i + \textnormal{j} Q_i$, complex coupling admittance $\underline{Y}_{\textnormal{coupl}} = \frac{1}{R_{\textnormal{S}}+\textnormal{j} \omega L_{\textnormal{S}}}$ (or, for the inverters $\underline{Y}_{\textnormal{coupl}} = \frac{1}{\textnormal{j} \omega L_{\textnormal{C}}}$) and complex voltage $\underline{V_i} = V_i \textnormal{e}^{\textnormal{j} \Delta \theta_i}$, we obtain the complex voltage $\underline{V}_{\textnormal{grid},i} = V_{\textnormal{grid},i} \textnormal{e}^{\textnormal{j} \Delta \theta_{\textnormal{grid},i}}$ between VISMA stator or inverter filters and grid as $
\underline{V}_{\textnormal{grid},i}(t) = \frac{\left| \underline{V}_i(t) \right|^2}{\underline{\overline{V}}_i(t)} - \frac{\underline{\overline{S}}_i(t)}{3 \underline{\overline{Y}}_{\textnormal{coupl}} \underline{\overline{V}}_i(t)}
$.

\section{Problem Statement}
\subsection{Optimisation Constraints \label{opt_constr}}
\noindent The objective of the optimisation is to find parameters $J$, $k_\T{d}$, $T_\T{d}$, and $K_\T{I}$ for the VISMA which positively influence the overall system behaviour after a perturbation. In order to avoid undesired or physically impossible behaviour, the optimisation variables have to be bounded by user-defined constraints. 

The first constraint assures that the VISMA does not react faster than the other inverters.
Therefore, we investigate the dynamics of the VISMA (cf., \eqref{eq:VISMAdomega2}-\eqref{eq:VISMAdamping}). For the purpose of deriving a simple model as reference for the optimisation constraints, the following simplifications are used:
Only the machine model is investigated, grid and stator equations, differential equations of voltage dynamics and secondary control are not considered.
Taking $P_1(\boldsymbol{V}(t), \Delta \boldsymbol{\theta}(t))$ as system input $u(t)$, linearising around the equilibrium point $u^{\ast} = P_{\textnormal{nom},1}$, $\omega_1^{\ast} = \omega_{\textnormal{nom}}$ and $d^{\ast} = - \omega_{\textnormal{nom}}$, and applying the Laplace transform gives the transfer function:
\begin{equation}
G_{\mathrm{VISMA, lin}}(s) = - \frac{k_{\mathrm{P},1} (T_{\mathrm{d}} s + 1)}{\frac{1}{\Omega^2} s^2 + \frac{D}{2 \Omega} s + 1}, \hspace{0.5cm} 
c := \frac{1}{k_{\mathrm{P},1} \omega_{\mathrm{nom}}}, \hspace{0.5cm} D := \frac{\frac{1}{c} (k_{\mathrm{d}}+J) +  T_{\mathrm{d}}}{2 \sqrt{\frac{1}{c} J T_{\mathrm{d}}}}, \hspace{0.5cm}
\Omega := \frac{1}{\sqrt{\frac{1}{c} J T_{\mathrm{d}}}}
\label{eq:linearizedSimpleVISMAtf} 
\end{equation}
The poles are real because $D > 1$ for any choice of parameters (proof omitted),  therefore:
\begin{equation}
s_{\textnormal{pole},1} = - \Omega \left( D + \sqrt{D^2 - 1} \right), \hspace{0.5cm} 
s_{\textnormal{pole},2} = - \Omega \left( D - \sqrt{D^2 - 1} \right) \label{eq:linearizedSimpleVISMApoles}
\end{equation}

From linear system theory it is known that $\tau_{1/2} = -\frac{1}{s_{\textnormal{pole},1/2}}$ determines the exponential decay rate.
This results in the constraint $\max_{i} \left( T_i \right) \leq \min \left( \tau_1, \tau_2 \right)$, where $T_i$ is the time constant of the regular inverters. Assuming stable system configurations, i.e., $s_{\textnormal{pole},1}, s_{\textnormal{pole},2} < 0$, we conclude that $\tau_1 < \tau_2$, and  therefore:
\begin{equation}
\max_{i} \left( T_i \right) \leq \frac{1}{\Omega \left( D + \sqrt{D^2 - 1} \right)} \label{n_Ti}
\end{equation}

A second constraint defines an upper bound for the parameter $K_\mathrm{I}$ of the integral controller \eqref{eq:integral_controller} based on the design preference that integral control action should occur only after the first part of the transient caused by droop control is finished.
For the simplified model \eqref{eq:linearizedSimpleVISMAtf} more than 95 \% of the absolute value of the step size are reached after $3 \mathrm{max}_{i=1,2}(\tau_i)$ because of the exponential character of the linearised system's step response. The response time of the integral controller should be larger.
A lower bound for the response time $\frac{1}{K_{\mathrm{I}}}$ of the integral controller, by using that $x(t)$ is scaled by $\frac{1}{J \omega_{\mathrm{nom}}}$ in \eqref{eq:linearizedSimpleVISMAtf}, is given by:
\begin{equation}
K_{\mathrm{I}} \leq 
\frac{ J \omega_{\mathrm{nom}} }{3 \tau_2} = \frac{1}{3} \, J \, \omega_{\mathrm{nom}} \, \Omega \left( D - \sqrt{D^2 - 1} \right) \label{eq:KI}
\end{equation} 

\subsection{Cost Functional \label{sec:Problem_statement_Analystical_investigation}}
\noindent The cost functional to be optimised contains three parts with parameters $\alpha > 0$ and $\beta > 0$:  
\begin{equation}
E [\Delta f, \Delta V, \delta_\T{f}, \delta_\T{V}, \alpha, \beta] =  t_{\mathrm{final}} + \alpha \cdot (k_{\mathrm{d}} + J) + 
\underbrace{(\Delta f / \delta_\T{f} + \Delta V / \delta_\T{V})}_{=: \Sigma}/\beta
\label{E_peaks}
\end{equation}

$\alpha$ and $\beta$ allow to shift the focus of the optimisation. First, we want the transient after a perturbation to be as short as possible, i.e., 
$t_{\mathrm{relax}} \rightarrow \mathrm{min}$. The time $t_\T{relax}$ is the relaxation time of the system after a perturbation. It is defined as the largest time of the moments, when the frequencies reach 49.999 Hz again. 
The time interval $t_\T{final} = t_\T{relax} - t_0$ considers the moment $t_0$ when the jump in load occurs.
Second, we consider the peak depth in the transients of frequency and voltage. We want them to be as small as possible to avoid damage on electronic devices, i.e., $\Delta f / \delta_\T{f} + \Delta V / \delta_\T{V} \rightarrow \min$, where $\Delta f = \max_{ \{i \in \{1,2,3\} \} ,\ t > t_0\}} |f_i(t) - f_i(t_0)|$ is the maximum frequency deviation, 
$f_i$ being the frequency at node~$i$. The max.~voltage deviation is $ \Delta V = \max_{\{i \in \{1,2,3,4\},\ t > t_0\}} |V_{\textnormal{grid},i}(t) - V_{\textnormal{grid},i}(t_0)|$, 
with grid voltage $V_{\textnormal{grid},i}(t)$. 

Third, a trade-off must be found between the required storage capacity of the VISMA, which should be as small as possible, and the energy that is used to keep up its virtual inertia.
By setting $M_\T{mech}=0$ and integrating \eqref{eq:VISMAdomega2} over time, the energy provided to or taken from the microgrid by the VISMA is:
\begin{equation}
E_\T{VISMA} = - \frac{1}{2} \left(J+k_{\mathrm{d}}\right) \left(\omega^2(t_2) -\omega^2(t_1) \right) +T_{\mathrm{d}}\ \int_{t_1}^{t_2} \omega_i(t) \ \dot{M}_{\mathrm{d}}(t) \,\mathrm{d}t
\end{equation}

$T_{\mathrm{d}}$ is responsible for scaling the energy loss due to damping and the other part of the energy is dominated by $J+k_{\mathrm{d}}$.
To avoid unnecessary large storage capacities $k_{\mathrm{d}} + J$ is minimised.
The constraints of the optimisation problem allow some insights into the structure of the optima in advance.
Choosing a very small $k_{\mathrm{d}}$ and $T_{\mathrm{d}}$ close to $\max_{i=1,2}{T_i}$ gives values as close as possible to the lower bound of \eqref{n_Ti}.
On the one hand this indicates that an optimisation with focus on keeping the transient behaviour of the VISMA close to those of regular inverters (i.e., minimising the virtual inertia $J+k_{\mathrm{d}}$ and the transient time $t_{\mathrm{final}}$) will give results close to $T_{\mathrm{d}} \approx \max_{i=1,2}{T_i}$, $k_{\mathrm{d}} \approx 0$ and $J \approx c \max_{i=1,2}{T_i}$.
Minimising mainly $t_{\mathrm{final}}$ on the other hand leads to a large value within limits given by \eqref{eq:KI}, namely $K_{\mathrm{I}} \approx \frac{1}{3 k_{\mathrm{P},1}}$.

\section{Implementation}
\noindent The energy landscape (see \FG{energy_landscape}) is rough with many local minima, in particular due to a stochastic term needed for the initial conditions when solving the differential equations, preventing the application of standard, e.g., gradient-based, methods. Instead, we use Parallel Tempering here, which also works for harder optimization problems, but is easy to implement.
\begin{figure}[ht]
	\centering
	\includegraphics[width=0.4\textwidth]{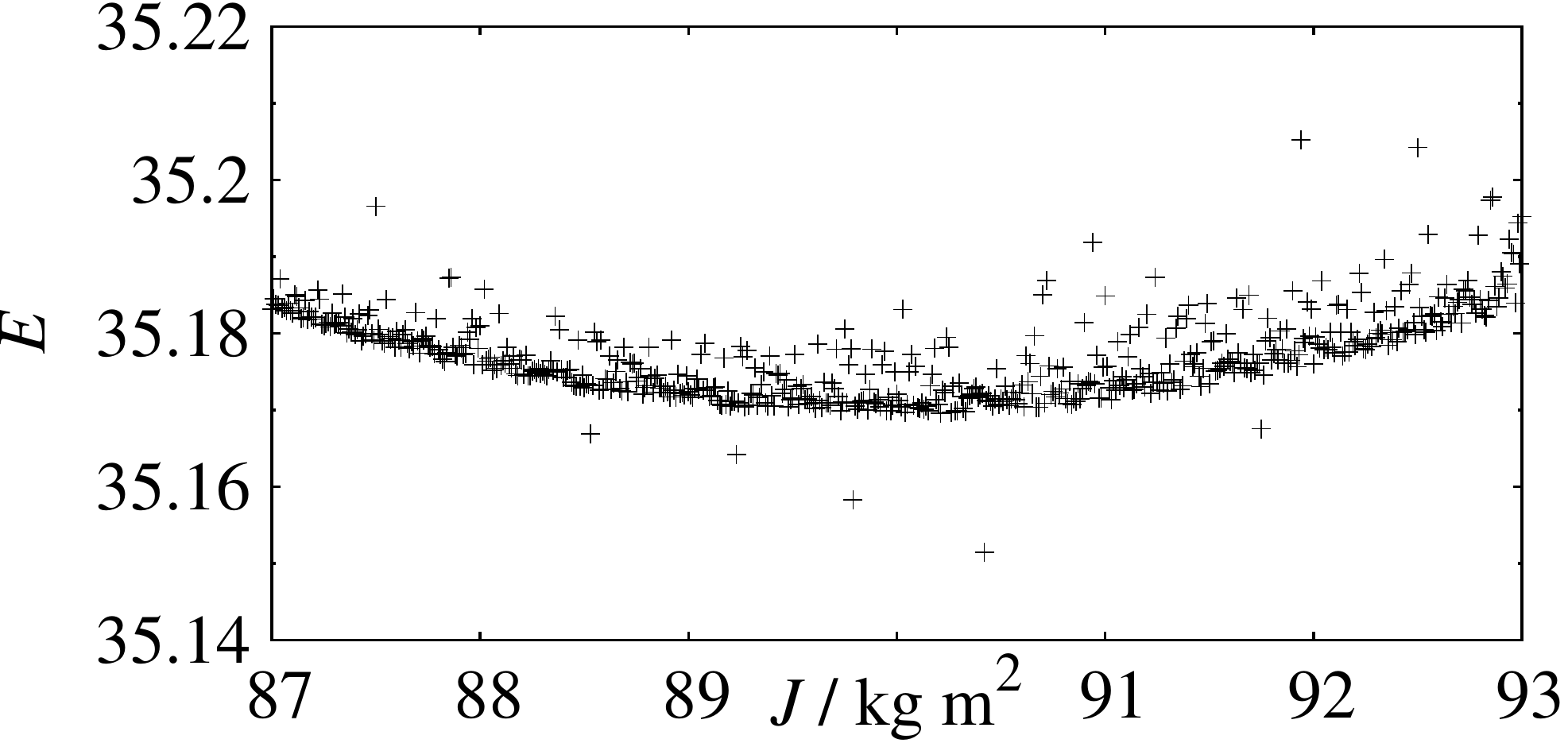}
	\caption{2D projection of the energy landscape (value of cost functional $E$ \eqref{E_peaks}) for scenario 1, close to Min.~$\# 2$. 
		The parameters $T_\T{d} = 0.6$, $k_\T{d} = 2.6 \cdot 10^{-4}$ , and $K_\T{I} = 1060$ are fixed, 
		whereas $J$ is varied. \label{energy_landscape}}
\end{figure}
\vspace*{-6mm}

\subsection{Parallel Tempering \label{Parallel_Tempering}}
\noindent The optimisation algorithm works as follows.
The configurations of the system are sampled according to the Boltzmann probability distribution $P(E_i) = 1/Z \, \exp(-E_i/\Theta)$, where $Z$ is a normalization constant, $E_i$ is the energy of configuration~$i$ and $\Theta$ an artificial temperature.
This is achieved via a special Monte Carlo~(MC) sampling, the Metropolis alogrithm \cite{Metropolis53}, where in each iteration, a new candidate configuration with corresponding energy $E_2$ is created and accepted with probability:
\begin{equation}
 p_\T{Metr} = \min \left\{ 1, e^{-(E_2 - E_1)/\Theta} \right\}
 \label{p_Metrop}
\end{equation}

$E_1$ is the energy of the current state. As \eqref{p_Metrop} fulfills \textit{detailed balance}, sampling according to a Boltzmann distribution is ensured.
From physics we know that for $\Theta \rightarrow 0$, the energy obtains a minimum $E \rightarrow E_\T{min}$. This leads to the idea of \textit{Simulated Annealing} \cite{Kirkpatrick83}, where the temperature of an MC simulation is gradually decreased until a minimum is found. This approach can get stuck in a local minimum. An improvement is to simulate the system at various temperatures $\Theta_i$. This can be done by Parallel Tempering \cite{Swendsen86,Hukushima96}, where a random walk in temperature space is performed. To preserve detailed balance and equilibrium for an infinite number of iterations, the Metropolis criterion \cite{Hukushima96} with energies $E(\cdotp)$ is used:
\begin{equation}
 p_\T{Swap} = \min \left\{ 1, \exp \left(\left[\frac{1}{\Theta_k} - \frac{1}{\Theta_{k+1}} \right] [E(y_k)-E(y_{k+1})] \right) \right\}
 \label{p_Swap}
\end{equation}

Two neighbouring configurations with temperatures $\Theta_k, \Theta_{k+1} \; (k \in [1,n-1])$ can be exchanged.
In each such \textit{swap}, $k \in \{1,2,3, \ldots , n-1\}$ is chosen at random with equal probability.

\subsection{Implementation of the Optimisation Algorithm \label{implementation} \label{Cost_func_sec}}
\noindent Within an optimisation procedure, one of two perturbation scenarios is considered. 
Both are based on a step in load. 
During the transient it is checked whether the usual frequency and voltage ranges are met (see \eqref{kp_kQ}). If these are not fulfilled, the parameter set is rejected, i.e., $E = \infty$. 
Before \eqref{E_peaks} is calculated in the simulation, the constraints \eqref{n_Ti} and \eqref{eq:KI} are checked. If they are not fulfilled, the parameter set is rejected.

We use from the GSL \cite{GSL}: A hybrid method (Newton and dogleg step) for solving the steady state equations, where the results are used as initial conditions for the Runge-Kutta-Fehlberg method to solve the differential equations. For parallelization we use OpenMPI \cite{OpenMPI}.
The 12 temperatures used in the simulations are $\Theta_i \in \{0.01,0.02,0.07,0.2,0.5,1,3,7,20,50,100,10^9\}$, where $10^9$ corresponds to the acceptance of every new state except the ones that violate \eqref{n_Ti}, \eqref{eq:KI} or lead to an unstable system state.
For each temperature $\Theta_i$ the MC sampling is performed in the following way:
\begin{enumerate}
 \itemsep0pt
 \item Calculate value of cost functional $E_1$ with given parameter set $\Phi = (J, k_\T{d}, T_\T{d}, K_\T{I})$.
 \item Choose one parameter $O$ of the four parameters from $\Phi$ with equal probability at random.
 \item Calculate $O' = O \cdot m$ with $m = |1 + R_\T{perc} \cdot r|$, where $r \in [-1,1]$ is a random number.
 \item Calculate new value of $E_2$ with modified parameter set $\Phi'$ in which $O$ has been replaced by $O'$.
 \item Accept the new parameter set $\Phi'$ with Metropolis probability \eqref{p_Metrop}.
\end{enumerate}

The steps 2.-5. are repeated $N_\T{params} \cdotp 2 = 4 \cdotp 2 = 8$ times.
After two such \textit{sweeps} are performed for each temperature, $n-1$ swap attempts are done. For each of these attempts the procedure is the following: 
First, choose a configuration $k \in [1,n-1]$ uniformly at random. Then, exchange the two configurations $y_k$ and $y_{k+1}$ with the swap probability given by \eqref{p_Swap}.
For each parameter set $(\alpha, \beta)$, 200 swaps are performed with $R_\T{perc} = 0.8$. The minimum of $E$ is found by taking the minimum of all temperatures $\Theta_i$ leading to $\Phi_\T{min}$.
Another simulation \cite{Hartmann15} with $\Phi_\T{min}$ as initial parameter set is started with 200 swaps and $R_\T{perc} = 0.4$.

\section{Results}
\subsection{Optimisation Results for Different Perturbation Scenarios}
\noindent For the optimisation a microgrid in a radial topology and in island mode is considered, 
see Fig.~\ref{fig:Disturbance_scenario_setup}. We neglect ohmic grid losses and focus on the frequency peak, i.e., $\delta_\T{V} = 10^{40}$. The perturbation is a jump in active load. 
Table~\ref{tab:System_parameters_fixed} shows the parameters that remain the same for all scenarios.

\begin{table}[ht]
	\centering
	\caption{Parameters Used for the Optimisation in All Scenarios. Values for $T_{2/3}$ taken from \cite{DroopControlSchiffer2014}
		\label{tab:System_parameters_fixed}}
	\begin{tabular}{c|c|c|c|c|c|c|c|c|c}
		$L_{\textnormal{lines}}$ & $L_{\textnormal{S}}$ & $R_{\textnormal{lines}}$ & $R_{\textnormal{S}}$ & 
		$L_{\textnormal{C}}$ & $T_1$ & $T_{2/3}$ & $k_{\mathrm{V}}$ & $Q_{\textnormal{nom},1/2/3}$ & $K_{\mathrm{awu}}$ \\
		\hline		
		1.514 mH & 42.0 mH & 0.0 $\Omega$ & 0.3 $\Omega$ &
		1.8 mH & 0.01 s & 0.5 s & 10.0  $\frac{\mathrm{V}}{\mathrm{V}}$ & 
		0.0 var& 1.0 $\frac{1}{\mathrm{s}}$ 
	\end{tabular}
\end{table}

\subsubsection{Scenario 1: Symmetric nominal power, load jump of 3 kW \label{scen1}}
\noindent In this scenario, a step decrease of the load power is done. See Table \ref{tab:System_parameters_scenario1} for all parameters. 
In the last three columns of Table \ref{tab:ResultsScenario1_f_startMin} the values of the parts of the cost functional~\eqref{E_peaks} are given, which reflect on which part of the functional the optimisation has been focused.
For the first minimum ($\#1$) the three parts of \eqref{E_peaks} have the same weight, the second minimum (min.) focuses on $t_\T{final}$, the third on the value of $J+k_\T{d}$, and the fourth on small frequency peaks ($\Sigma$).
For the second, third and forth min.~the focus is reflected in the result, i.e., min.~$\# 2$ has the smallest value of $t_\T{final}$, min.~$\# 3$ the smallest of $J+k_\T{d}$. Min.~$\# 4$ gives a comparably small value for $\Sigma$, but apparently only a local min.~was found, since in min.~$\# 2$ it is even smaller. 

\begin{table}[hb]
	\centering
	\caption{Parameters Used for the Optimisation in Scenario 1 \label{tab:System_parameters_scenario1}}
	\begin{tabular}{c|c|c|c|c|c}
		$S_{1/2/3}$ & $P_{\textnormal{nom},1/2/3}$ & $k_{\textnormal{P},1/2/3}$ & $k_{\textnormal{Q},2/3}$ &  $P_{\textnormal{load}}$ before jump & $P_{\textnormal{load}}$ after jump \\
		\hline 4000.0 VA &  500.0 W & 3.1416$\cdot 10^{-4}$ $\frac{\textnormal{rad}}{\textnormal{s VA}}$ & 5.75$\cdot 10^{-3}$ $\frac{\textnormal{V}}{\textnormal{VA}}$ &  1500.0 W & 4500.0 W \\ 
	\end{tabular}
\end{table}

\begin{table}[ht]
	\centering
	\caption{Optimal Parameter Sets for Scenario 1 ($R_\T{perc} = 0.4$, Focus on Frequency Peak: 
	$\delta_\T{f} = 0.05$, $\delta_\T{V} = 10^{40}$).
	Errors of $E$ Resulting from 50 Runs with Different Initial Conditions. \label{tab:ResultsScenario1_f_startMin}}
	\resizebox{\textwidth}{!}{\begin{tabular}{r|r|r|r|r|r|r|r|r|r|r|r|r}
	\# & $J$ & $k_{\mathrm{d}} / 10^{-4}$ & $T_{\mathrm{d}}$ & $K_{\mathrm{I}}$ & $E$ & $\alpha$ & $\beta$ & 
	$J + k_{\mathrm{d}}$ & $\Sigma$ & $t_{\mathrm{final}}$ & $\alpha (J + k_{\mathrm{d}})$ & $\Sigma/\beta$ \\
		\hline
1 & 5.0895 & 1.1857 & 0.5029 & 1054.56 & 108.93(6) & 7 & 0.027 & 5.090 & 0.994 & 36.483 & 35.627 & 36.820 \\
2 & 91.479 & 2.5800 & 0.5917 & 1060.97 & 35.12(2) & 0.07 & 2.7 & 91.480 & 0.817 & 28.415 & 6.4036 & 0.3026 \\
3 & 5.0692 & 1.0071 & 0.5163 & 975.67 & 3624.89(9) & 700 & 0.027 & 5.069 & 1.000 & 39.379 & 3548.494 & 37.026 \\
4 & 50.894 & 10.1498 & 1.2539 & 1053.54 & 3425(46) & 7 & $2.7 \cdot 10^{-4}$ & 50.895 & 0.820 & 32.913 & 356.265 & 3036.54 \\
	\end{tabular}}
\end{table}

In Sec.~\ref{sec:Problem_statement_Analystical_investigation} it was stated that $T_{\mathrm{d}} \approx \max_{i=1,2}{T_i} = 0.5$, $k_{\mathrm{d}} \approx 10^{-4}$ (since it is bounded by this value) and $J \approx c \max_{i=1,2}{T_i} \approx 10.13 \cdot 0.5 \approx 5.07$ in cases where $J+k_{\mathrm{d}}$ and $t_{\mathrm{final}}$ are minimised with equal weighting, and $K_{\mathrm{I}} \approx \frac{1}{3 k_{\mathrm{p},1}} \approx 1061.03$, if the minimisation focuses on $t_{\mathrm{final}}$. These values are close to min.~$\# 1$ and $\# 3$, and min.~$\# 2$ for focussing on $t_{\mathrm{final}}$.
Weightings on other parts of the cost functional (e.g., min.~$\# 4$) lead to minima which are further away from the bounds of the optimisation constraints. 

\begin{figure}[ht]
 \centering
 \subfigure[Frequencies.]{\includegraphics[width=0.37\textwidth]{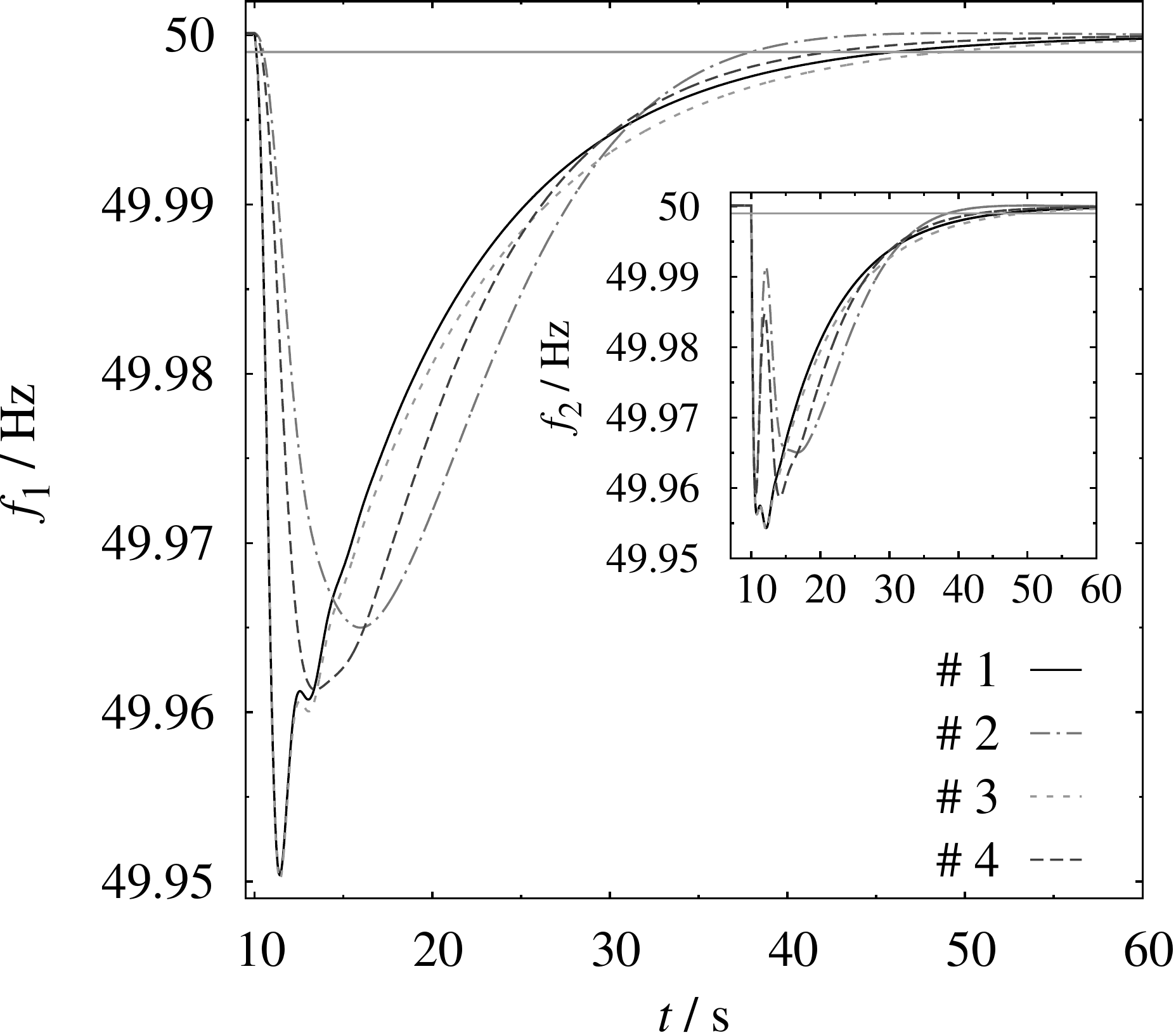}}
 \hspace{6mm}
 \subfigure[Voltages.]{\includegraphics[width=0.37\textwidth]{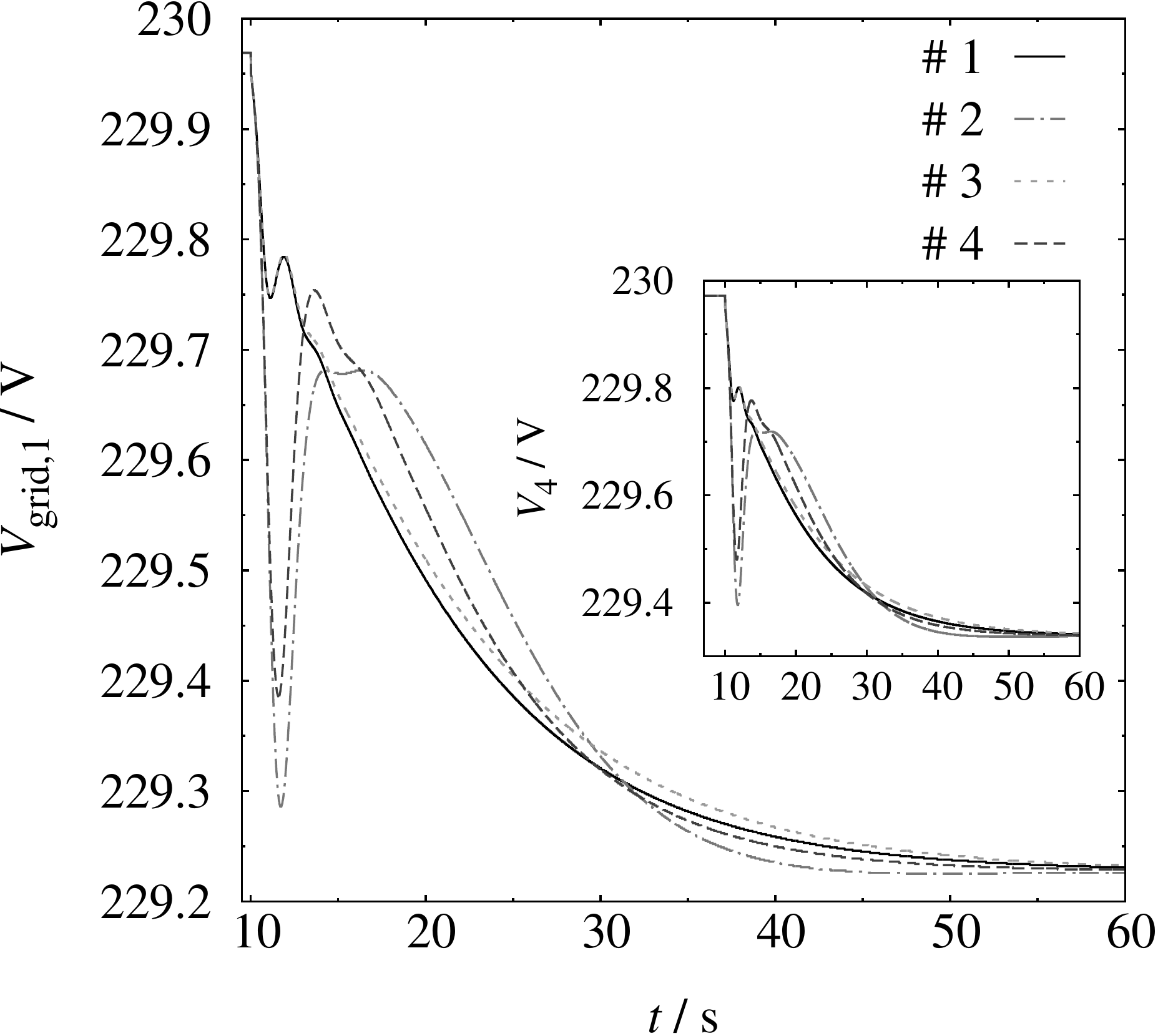}}
 \caption{Comparison of the four minima for scenario 1 given in Table \ref{tab:ResultsScenario1_f_startMin}. \label{fig:Comparison_of_Min_Scenario1}}
\end{figure}

The results from Table \ref{tab:ResultsScenario1_f_startMin} are visualized (see Fig.~\ref{fig:Comparison_of_Min_Scenario1}) by a comparison for VISMA frequencies and voltages for all four minima.
The remaining grid values are shown only for the first minimum, see Fig.~\ref{fig:Scenario1Min1SystemBehaviour}.  

\subsubsection{Scenario 2: Different nominal powers and load jump of 7 kW}
\begin{figure}[ht]
	\centering
	\begin{minipage}{0.48\textwidth}
	 \includegraphics[width=0.85\textwidth]{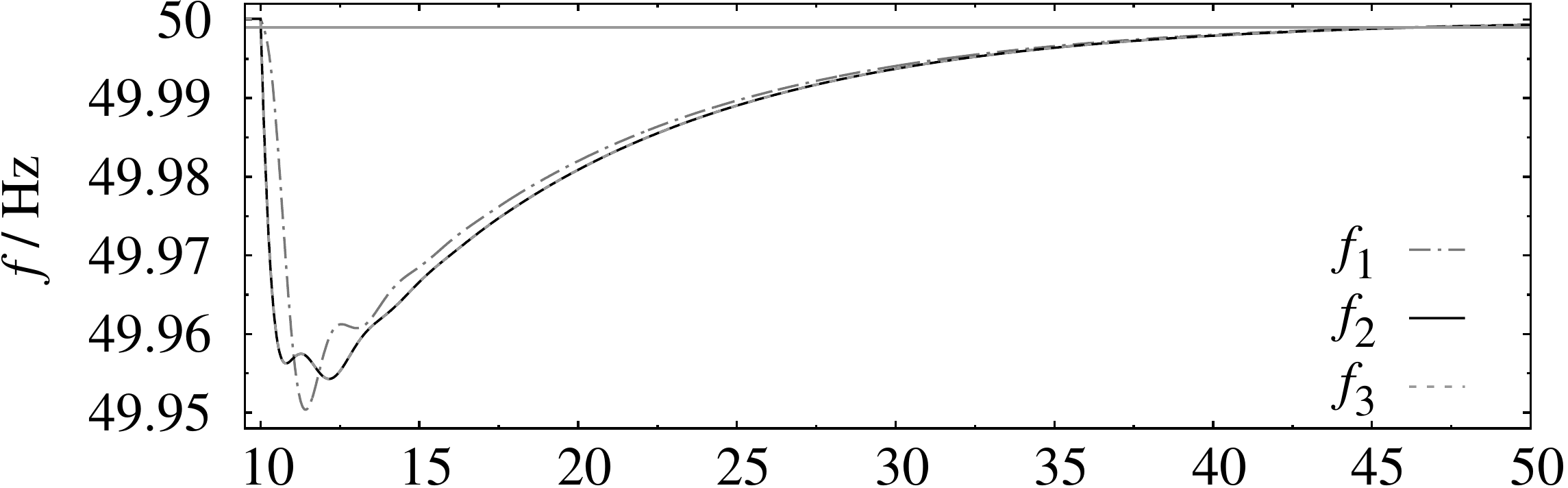}
	 \includegraphics[width=0.85\textwidth]{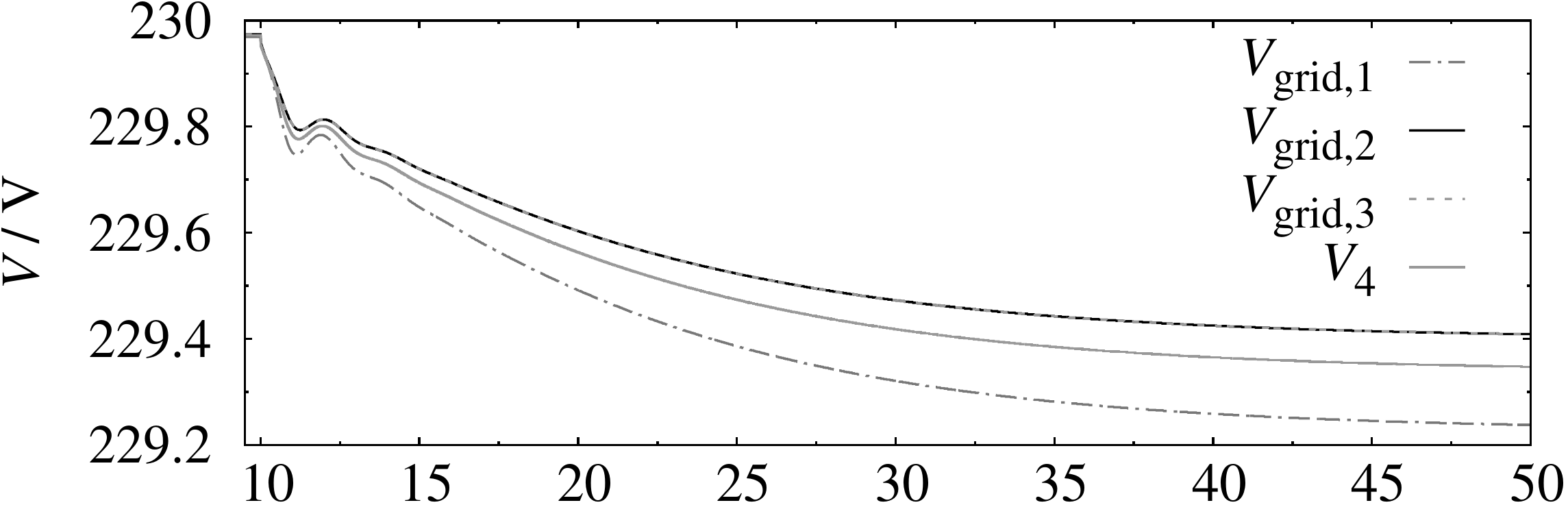}
	 \includegraphics[width=0.85\textwidth]{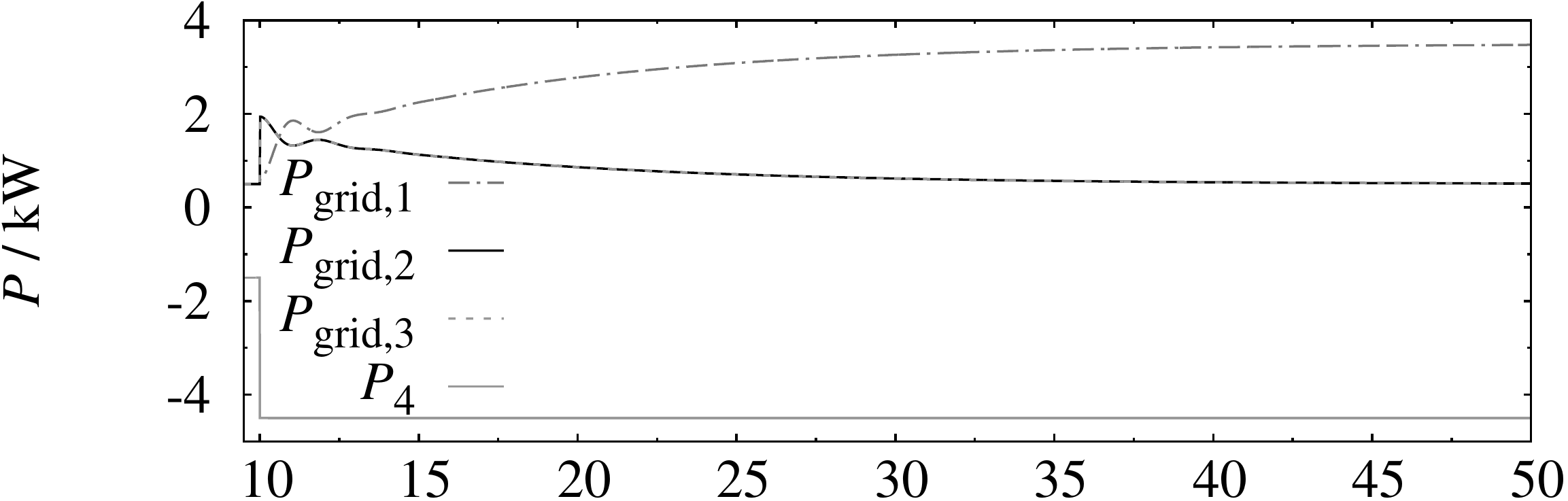}
	 \includegraphics[width=0.85\textwidth]{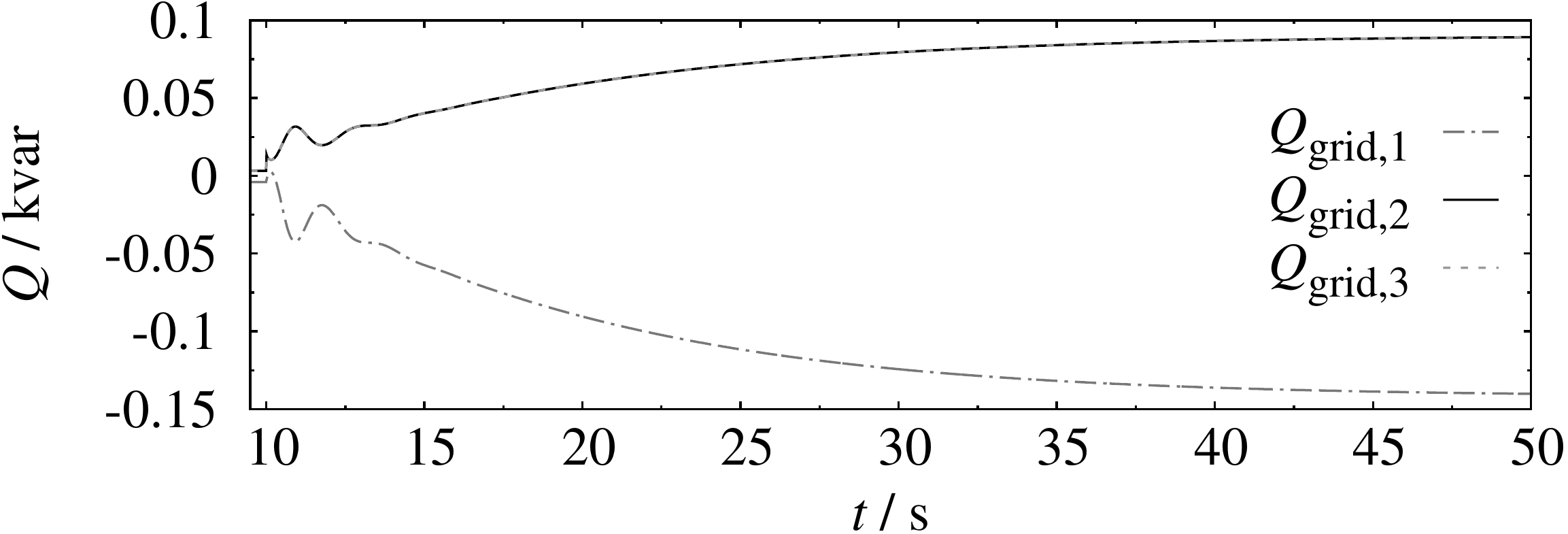}
	 \captionof{figure}{Minimum \#1 of scenario 1. \label{fig:Scenario1Min1SystemBehaviour}}
	\end{minipage}
	\hspace{2mm}
	\begin{minipage}{0.48\textwidth}
	\includegraphics[width=0.85\textwidth]{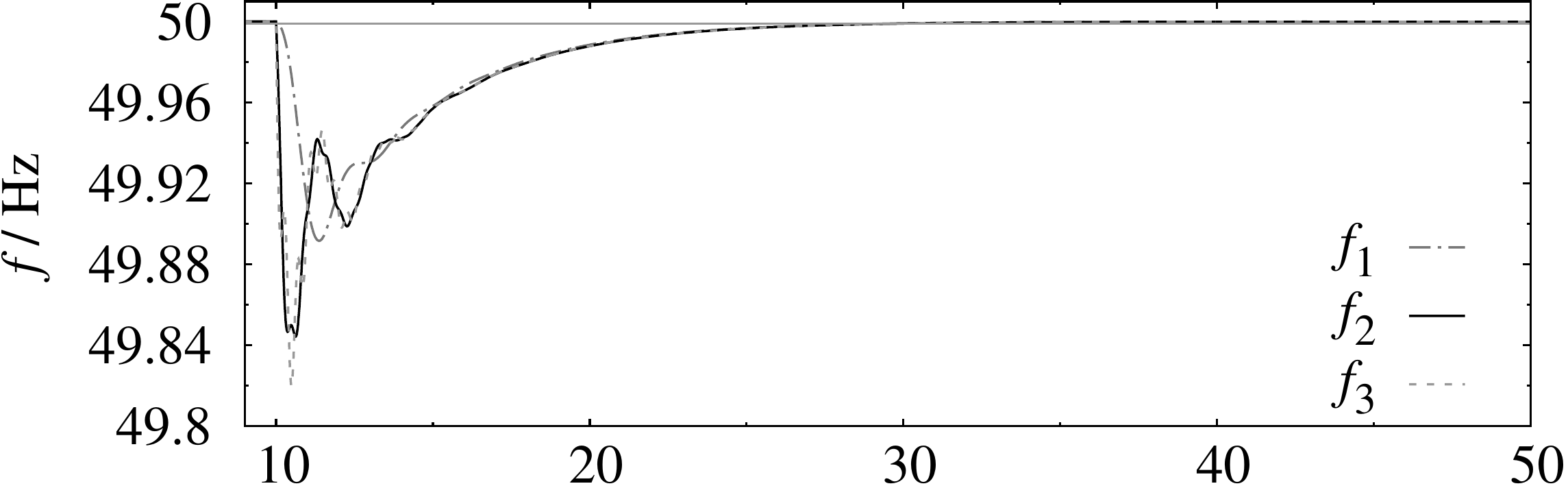}
	\includegraphics[width=0.85\textwidth]{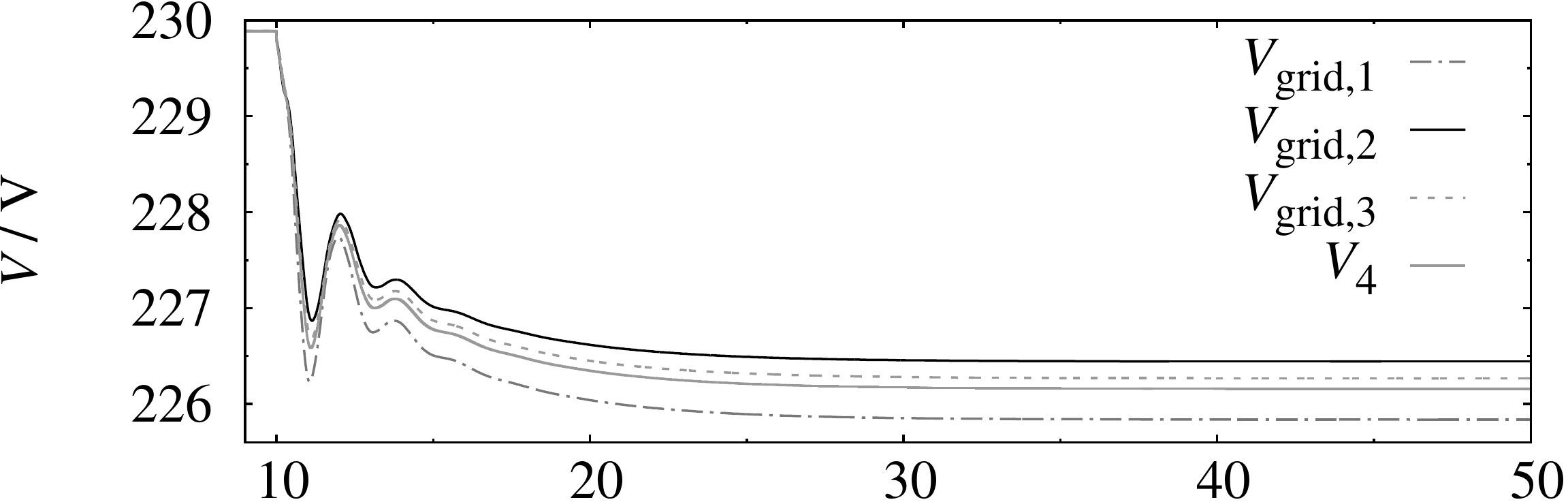}
	\includegraphics[width=0.85\textwidth]{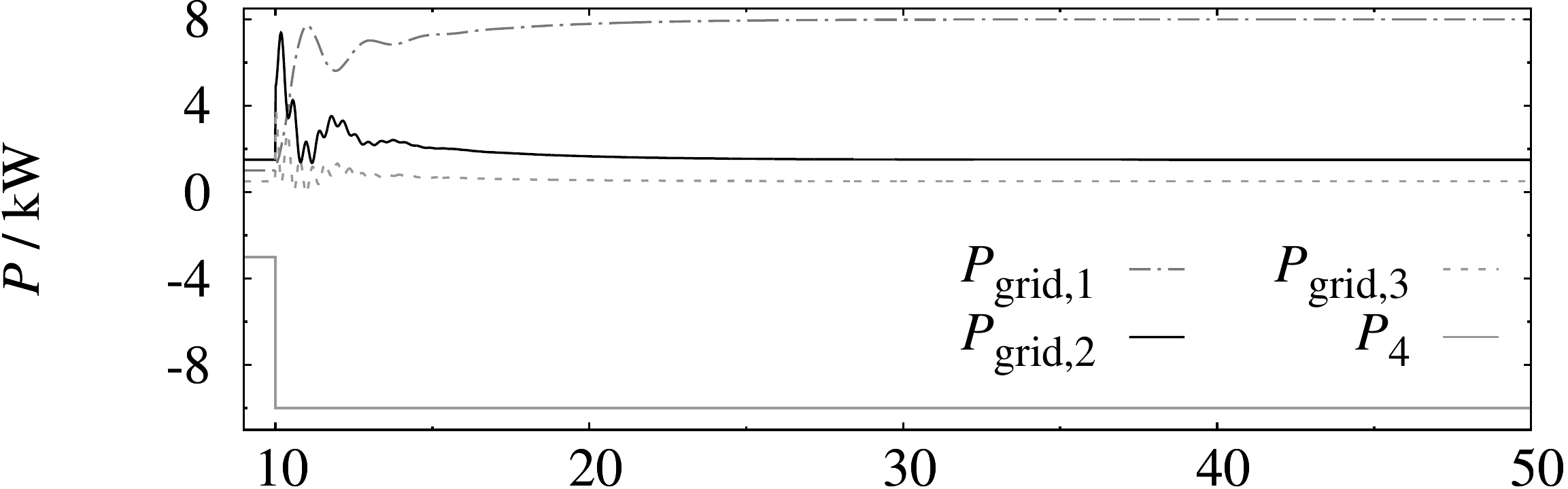}
	\includegraphics[width=0.85\textwidth]{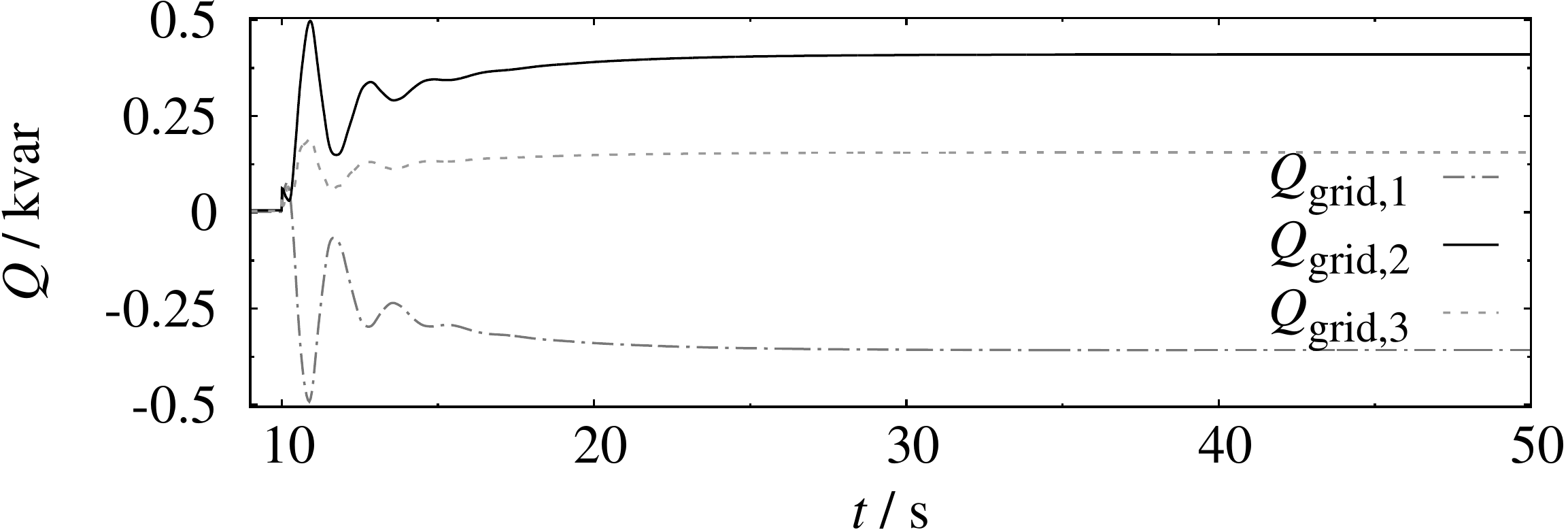}
	\captionof{figure}{Minimum \#1 of scenario 2. \label{fig:Scenario3Min1SystemBehaviour}}
	\end{minipage}
\end{figure}

\noindent In this scenario, we assume different rated and nominal active powers of the devices (see Table \ref{tab:System_parameters_scenario3}).
The optimal parameter sets obtained for this setup are shown in Table \ref{tab:ResultsScenario3_startMin}.
Fig.~\ref{fig:Scenario3Min1SystemBehaviour} shows the system behaviour for min.~$\# 1$ of Scenario 2. Different time constants of the VISMA and the two inverters cause a different system behaviour than in Scenario~1. 
The VISMA's reaction is very slow due to its high virtual inertia (for min.~$\#1$, $J \approx 11.5$). Directly after the load jump the inverters have to balance the sudden power demand. This forces them to provide active power at a value ($\sim 8$ kW) above their nominal rated power values $S_{2/3}$. Inverters for island grids allow this for a short amount of time. 

\begin{table}[htb]
	\centering
	\caption{Parameters Used for the Optimisation in Scenario 2 \label{tab:System_parameters_scenario3}}
	\begin{tabular}{c|c|c|c|c|c|c|c}
		$S_{1}$ & $S_{2}$ & $S_{3}$ & $P_{\textnormal{nom},1}$ & $P_{\textnormal{nom},2}$ & $P_{\textnormal{nom},3}$ & $k_{\textnormal{P},1}$ & $k_{\textnormal{P},2}$ \\
		\hline
		9.0 kVA & 3.0 kVA & 1.0 kVA & 1.0 kW & 1.5 kW & 0.5 kW & 1.3963$\cdot 10^{-4}$ $\frac{\textnormal{rad}}{\textnormal{s VA}}$ & 		4.1888$\cdot 10^{-4}$ $\frac{\textnormal{rad}}{\textnormal{s VA}}$ \\
	\end{tabular}
	\vspace*{2mm}
	
	\begin{tabular}{c|c|c|c|c}
		$k_{\textnormal{P},3}$ & $k_{\textnormal{Q},2}$ & $k_{\textnormal{Q},3}$ &  $P_{\textnormal{load}}$ before jump& $P_{\textnormal{load}}$ after jump \\
		\hline 12.5664$\cdot 10^{-4}$ $\frac{\textnormal{rad}}{\textnormal{s VA}}$ & 7.67$\cdot 10^{-3}$ $\frac{\textnormal{V}}{\textnormal{VA}}$ & 23.0$\cdot 10^{-3}$ $\frac{\textnormal{V}}{\textnormal{VA}}$ & 3.0 kW & 10.0 kW \\
	\end{tabular}
\end{table}

\begin{table}[!h]
	\centering
	\caption{Optimal Parameter Sets for Scenario 2 ($R_\T{perc} = 0.4$, Focus on Frequency Peak: 
	$\delta_\T{f} = 0.2$, $\delta_\T{V} = 10^{40}$).
	Errors of $E$ Resulting from 50 Runs with Different Initial Conditions. \label{tab:ResultsScenario3_startMin}}
	\resizebox{\textwidth}{!}{\begin{tabular}{r|r|r|r|r|r|r|r|r|r|r|r|r}
	\# & $J$ & $k_{\mathrm{d}} / 10^{-4}$ & $T_{\mathrm{d}}$ & $K_{\mathrm{I}}$ & $E$ & $\alpha$ & $\beta$ & 
	$J + k_{\mathrm{d}}$ & $\Sigma$ & $t_{\mathrm{final}}$ & $\alpha (J + k_{\mathrm{d}})$ & $\Sigma/\beta$ \\
		\hline
1 & 11.4986 & 1.1595 & 0.5035 & 2379.26 & 59.26(1) & 1.7 & 0.045 & 11.499 & 0.902 & 19.671 & 19.548 & 20.046 \\
2 & 59.3043 & 4.2760 & 0.9356 & 2387.04 & 14.99(1) & 0.017 & 4.5 & 59.305 & 0.887 & 13.781 & 1.0082 & 0.1971 \\
3 & 11.4076 & 1.0139 & 0.5064 & 2348.85 & 1979.32(1) & 170 & 0.045 & 11.408 & 0.902 & 19.967 & 1939.309 & 20.040 \\
4 & 16.8974 & 13.2759 & 0.6524 & 2382.76 & 2039.5(5) & 1.7 & $4.5 \cdot 10^{-4}$ & 16.899 & 0.896 & 19.076 & 28.728 & 1991.7 \\
	\end{tabular}}
\end{table}

Despite the different system behaviour compared to Scenario~1, the effect of the optimisation with different weights is clearly reflected in the values of the cost functional (see Table \ref{tab:ResultsScenario3_startMin}). 
The very small $t_{\mathrm{final}}$ for min.~$\# 2$ is due to an overshoot in frequencies (figure not shown): After reaching $49.999$~Hz for the first time, $f_1$ leaves the range of $\left[49.999, 50.001 \right]$ again before returning to $50$~Hz.
Optimisations for the same scenarios with ohmic losses in lines have been performed (results not shown). Though transients show different characteristics, almost the same parameter sets are found.

\section{Conclusion}
\noindent Parameters of a VISMA in an islanded microgrid with radial topology containing two inverters and a load have been optimised using Parallel Tempering. 
By varying additional parameters in the cost functional the focus of the optimisation was shifted. For two perturbation scenarios minima of the cost functional were found which are stable solutions within the prescribed boundaries.
The results show that this optimisation procedure is in general applicable to the task of parameter optimisation in a microgrid. 
It is also shown that through the proper setting of the VISMA's parameters its functionality can be adapted to different participants in the grid. The values obtained by analytical investigation in \ref{sec:Problem_statement_Analystical_investigation} seem to offer a "rule of thumb" for a good parameter region. However, the effects of each parameter in complex situations are not obvious. For other topologies, devices or disturbances, totally different parameter sets might be needed. In order to find a good parameter set for a given microgrid setup, various disturbances should be analyzed. 

In future research, other forms of the constraints for the VISMA parameters and cost functionals should be investigated. Larger microgrids with other power generating systems can be put under scrutiny. The proposed approach could be transferred to the optimisation of parameters in other control strategies.
An extension would be the study of different disturbance scenarios, where one uses e.g., a series of steps in load. Finally, theoretical results should be confirmed by measurements in an appropriate laboratory.

\section*{Acknowledgment}
\noindent Simulations were performed at the HERO cluster of the Carl von Ossietzky Universit\"at Oldenburg, funded by the DFG through its Major Research Instrumentation Programme (INST 184/108–1 FUGG) and the Ministry of Science and Culture of the Lower Saxony State. Financial support was obtained by the Lower Saxony research network "Smart Nord" which acknowledges the support of the Lower Saxony Ministry of Science and Culture through the "Nieders\"achsisches Vorab" grant programme (grant ZN2764).

\bibliography{Opt_MicroGrids_formatted}

\end{document}